
\magnification1200
\input amstex.tex
\documentstyle{amsppt}

\hsize=12.5cm
\vsize=18cm
\hoffset=1cm
\voffset=2cm

\def\DJ{\leavevmode\setbox0=\hbox{D}\kern0pt\rlap
 {\kern.04em\raise.188\ht0\hbox{-}}D}
\def\dj{\leavevmode
 \setbox0=\hbox{d}\kern0pt\rlap{\kern.215em\raise.46\ht0\hbox{-}}d}

\baselineskip=13pt
\def\hf{{\textstyle{1\over2}}}
\def\b{\beta}
\def\d{{\,\roman d}}

\def\f{\varphi}
\def\G{\Gamma}
\def\g{\gamma}

\def\s{\sigma}
\def\t{\theta}
\def\={\;=\;}

\def\R{\Re{\roman e}\,} 
\def\z{\zeta}

 \def\t{\theta}
\def\hf{{\textstyle{1\over2}}}

\def\f{\varphi}

\font\tenmsb=msbm10
\font\sevenmsb=msbm7
\font\fivemsb=msbm5
\newfam\msbfam
\textfont\msbfam=\tenmsb
\scriptfont\msbfam=\sevenmsb
\scriptscriptfont\msbfam=\fivemsb
\def\Bbb#1{{\fam\msbfam #1}}

\def \NN {\Bbb N}
\def \CC {\Bbb C}
\def \RR {\Bbb R}

\font\ff=cmr8

\baselineskip=13pt

\font\teneufm=eufm10
\font\seveneufm=eufm7
\font\fiveeufm=eufm5
\newfam\eufmfam
\textfont\eufmfam=\teneufm
\scriptfont\eufmfam=\seveneufm
\scriptscriptfont\eufmfam=\fiveeufm
\def\mathfrak#1{{\fam\eufmfam\relax#1}}

\font\tenmsb=msbm10
\font\sevenmsb=msbm7
\font\fivemsb=msbm5
\newfam\msbfam
     \textfont\msbfam=\tenmsb
      \scriptfont\msbfam=\sevenmsb
      \scriptscriptfont\msbfam=\fivemsb
\def\Bbb#1{{\fam\msbfam #1}}

\def \NN {\Bbb N}
\def \CC {\Bbb C}

\def \RR {\Bbb R}

  \def\rightheadline{{\hfil{\ff
  Some identities for the Riemann zeta-function II}\hfil\tenrm\folio}}

  \def\leftheadline{{\tenrm\folio\hfil{\ff
   Aleksandar Ivi\'c }\hfil}}
  \def\emptyheadline{\hfil}
  \headline{\ifnum\pageno=1 \emptyheadline\else
  \ifodd\pageno \rightheadline \else \leftheadline\fi\fi}

\font\ff=cmr8
\font\teneufm=eufm10
\font\seveneufm=eufm7
\font\fiveeufm=eufm5
\newfam\eufmfam
\textfont\eufmfam=\teneufm
\scriptfont\eufmfam=\seveneufm
\scriptscriptfont\eufmfam=\fiveeufm
\def\mathfrak#1{{\fam\eufmfam\relax#1}}

\font\tenmsb=msbm10
\font\sevenmsb=msbm7
\font\fivemsb=msbm5
\newfam\msbfam
\textfont\msbfam=\tenmsb
\scriptfont\msbfam=\sevenmsb
\scriptscriptfont\msbfam=\fivemsb
\def\Bbb#1{{\fam\msbfam #1}}

\def \NN {\Bbb N}
\def \CC {\Bbb C}

\def \RR {\Bbb R}

\def\b{\beta} 
 \def\d{\,{\roman d}}
\topmatter
\title
Some  identities for the Riemann zeta-function II
\endtitle
\author
Aleksandar Ivi\'c
\endauthor
\address
Katedra Matematike RGF-a Universiteta u Beogradu, \DJ u\v sina 7,
11000 Beograd, Serbia and Montenegro
\endaddress
\subjclass 11 M 06
\endsubjclass
\email {\tt aivic\@matf.bg.ac.yu, ivic\@rgf.bg.ac.yu}
\endemail
\dedicatory
Facta Universitatis (Ni\v s) {\bf20}(2005), 1-8.
\enddedicatory
\abstract
Several identities for the Riemann zeta-function
$\zeta(s)$ are proved. For example, if
$\f_1(x) := \{x\} = x - [x], \quad\f_n(x) := \int_0^\infty
\{u\}\f_{n-1}\bigl({x\over u}\bigr)\,{\d u\over u}\;(n\ge 2),$ then
$$
{\z^n(s)\over(-s)^{n}} = \int_0^\infty \f_n(x)x^{-1-s}\d x
\quad(s = \s + it,\;0 < \s < 1)
$$
and
$$
{1\over2\pi}\int_{-\infty}^\infty
{|\z(\s+it)|^{2n}\over(\s^2+t^2)^n}\d t = \int_0^\infty \f_n^2(x)x^{-1-2\s}\d x
\qquad(0<\s<1).
$$

\endabstract
\medskip
\endtopmatter

Let as usual $\z(s) = \sum_{n=1}^\infty n^{-s}\;(\R s > 1)$ denote
the Riemann zeta-function. This note is the continuation of the author's
work [6], where several identities involving $\z(s)$ were obtained.
The basic idea is to use properties of the Mellin transform
($f:[0,\infty)\rightarrow\RR$ )
$$
F(s) = \Cal{M}[f(x);s] := \int_0^\infty f(x)x^{s-1}\d x\qquad(s = \s+it,\,\s>0),
\leqno(1)
$$
in particular the analogue of the Parseval formula for Mellin
transforms, namely
$$
{1\over2\pi i}\int_{\s-i\infty}^{\s+i\infty}|F(s)|^2\d s
= \int_0^\infty f^2(x)x^{2\s-1}\d x.\leqno(2)
$$
For the conditions under which (2) holds, see e.g., [5] and [11]. If $\{x\}$
denotes the fractional part of $x$ ($\{x\} = x - [x]$, where $[x]$
is the greatest integer not exceeding $x$), we have the classical
formula (see e.g., eq. (2.1.5) of E.C. Titchmarsh [12])
$$
{\z(s)\over s} = - \int_0^\infty \{x\}x^{-1-s}\d x
= - \int_0^\infty \{1/x\}x^{s-1}\d x\quad(s = \s+it, \,0<\s<1).
\leqno(3)
$$
A quick proof is as follows. We have
$$\eqalign{
\z(s) &= \int_{1-0}^\infty x^{-s}\d[x] = s\int_{1}^\infty [x]x^{-s-1}\d x\cr&
= s\int_{1}^\infty ([x]- x)x^{-s-1}\d x + s\int_{1}^\infty x^{-s}\d x
= -s\int_{1}^\infty \{x\}x^{-s-1}\d x + {s\over s-1}.\cr}
$$
This holds initially for $\s>1$, but since the last integral is absolutely
convergent for $\s>0$, it holds in this region as well by analytic
continuation. Since
$$
s\int_{0}^1 \{x\}x^{-s-1}\d x = s\int_{0}^1 x^{-s}\d x = {s\over1-s}
\quad(0 < \s < 1),
$$
we obtain (3) on combining the preceding two formulae.
We note that (3) is a special case of the so-called M\"untz's formula
(with $f(x) = \chi_{[0,1]}(x)$, the characteristic function of the unit
interval)
$$
\z(s)F(s) = \int_0^\infty Pf(x)\cdot x^{s-1}\d x,\leqno(4)
$$
where the
M\"{u}ntz operator  $P$ is the linear operator defined formally
on functions $f:[0,\infty)\rightarrow\CC$ by
$$
Pf(x):=\sum_{n=1}^\infty f(nx) - \frac{1}{x}\int_0^\infty f(t)\d t.\leqno(5)
$$
Besides the original proof of (4) by M\"{u}ntz  [8], proofs are given by
E.C. Titchmarsh [12, Chapter 1, Section 2.11] and recently by L. B\'aez-Duarte [2].
The identity (4) is valid for $0<\s<1$ if $f'(x)$ is continuous, bounded in
any finite interval  and is $O(x^{-\b})$ for $x\to\infty$ where $\b>1$ is a
constant. The identity (3), which B\'aez-Duarte [2] calls the {\it
proto-M\"{u}ntz} identity, plays an important r\^ole in the approach to
the Riemann Hypothesis (RH, that all complex zeros of $\z(s)$ have real
parts equal to 1/2) via methods from functional analysis (see e.g., the works
[1]--[4] and [9]).

\medskip
Our first aim is to generalize (3). We introduce
the convolution functions $\f_n(x)$ by
$$
\f_1(x) := \{x\} = x - [x], \;\f_n(x) := \int_0^\infty \{u\}\f_{n-1}
\bigl({x\over u}\bigr)\,{\d u\over u}
\qquad(n\ge 2). \leqno(6)
$$
The asymptotic behaviour of the function $\f_n(x)$ is contained in

\bigskip
THEOREM 1. {\it If $n \ge 2$ is a fixed integer, then}
$$
\f_n(x) = {x\over(n-1)!}\log^{n-1}(1/x) + O\Bigl(x\log^{n-2}(1/x)\Bigr)
\qquad(0 < x < 1), \leqno(7)
$$
{\it and}
$$
\f_n(x) = O(\log^{n-1}(x+1)) \qquad (x\ge 1).\leqno(8)
$$

\medskip
{\bf Proof}. Using the properties of $\{x\}$, namely $\{x\} = x$
for $0 < x < 1$ and $\{x\} \le x$, one easily verifies (7) and (8)
when $n=2$. To prove the general case we use induction, supposing
that the theorem is true for some $n$. Then, when $0 < x < 1$,
$$
\f_{n+1}(x) = \int_0^x + \int_x^1 + \int_1^\infty = I_1 + I_2 + I_3,
$$
say. We have, by change of variable,
$$
\eqalign{
I_1 &= \int_0^x \{u\}\f_n\left({x\over u}\right){\d u\over u}
= \int_0^x \f_n\left({x\over u}\right)\d u
\cr& = x\int_1^\infty \f_n(v)\,{\d v\over v^2} = O(x).
\cr}
$$
By the induction hypothesis
$$
\eqalign{
I_2 &= \int_x^1 \{u\}\f_n\left({x\over u}\right){\d u\over u}
= \int_x^1 \f_n\left({x\over u}\right)\d u\cr&
= \int_x^1 \left\{{x\over(n-1)!u}\log^{n-1}\left({u\over x}\right)
+ O\left({x\over u}\log^{n-2}\left({u\over x}\right)\right)\right\}\d u\cr&
= {x\over(n-1)!}\int_1^{1/x}\log^{n-1}y{\d y\over y}
+ O\Bigl(x\log^{n-1}(1/x)\Bigr)\cr&
= {x\over n!}\log^n(1/x) + O(x\log^{n-1}(1/x)).\cr}
$$
Finally, since $\{x\} \le x$ and (8) holds,
$$
I_3 = \int_1^\infty \{u\}\f_n\left({x\over u}\right){\d u\over u}
\ll x\int_1^\infty \log^{n-1}\left({u\over x}\right){\d u\over u^2}
\ll x\log^{n-1}\left({1\over x}\right).
$$
The proof of (8) is on similar lines, when we write
$$
\f_{n+1}(x) = \int_0^1 + \int_1^x + \int_x^\infty = J_1 + J_2 + J_3
\qquad(x\ge 1),
$$
say, so that there is no need to repeat the details. By more elaborate
analysis (7) could be further sharpened.

\bigskip
THEOREM 2. {\it If $n \ge 1$ is a fixed integer, and $s = \s+it,\,0<\s<1$, then}
$$
{\z^n(s)\over(-s)^{n}} = \int_0^\infty \f_n(x)x^{-1-s}\d x.\leqno(9)
$$

\bigskip
Clearly (9) reduces to (3) when $n=1$. From Theorem 1 it transpires
that the integral in (9) is absolutely convergent for $0<\s<1$. By using (2)
(with $-s$ in place of $s$) we obtain the following

\bigskip
{\bf Corollary 1}. For $n\in\NN$ we have
$$
{1\over2\pi}\int_{-\infty}^\infty
{|\z(\s+it)|^{2n}\over(\s^2+t^2)^n}\d t = \int_0^\infty \f_n^2(x)x^{-1-2\s}\d x
\qquad(0<\s<1).\leqno(10)
$$

\bigskip
{\bf Proof of Theorem 2}. As already stated, (9) is true for $n=1$. The general
case is proved then by induction. Suppose that (9) is true for some $n$,
and consider
$$
{\z^{n+1}(s)\over(-s)^{n+1}} = \int_0^\infty \int_0^\infty
\{x\}\f_n(y)(xy)^{-1-s}\d x\,\d y\qquad(0<\s<1)
$$
as a double integral.
We make the change of variables $x = v, y = u/v$, noting that the absolute value
of the Jacobian of the transformation is $1/v$. The above integral becomes then
$$
\eqalign{&
\int_0^\infty \int_0^\infty \{v\}\f_n\left({u\over v}\right)u^{-1-s}v^{-1}\d u\d v\cr&
= \int_0^\infty \left(\int_0^\infty\{v\}\f_n\left({u\over v}\right){\d v\over v}\right)
u^{-1-s}\d u\cr&
= \int_0^\infty \f_{n+1}(x)x^{-1-s}\d x,\cr}
$$
as asserted. The change of integration is valid by absolute convergence, which is
guaranteed by Theorem 1.

\medskip
{\bf Remark 1}. L. B\'aez-Duarte kindly pointed out to me that the above
procedure leads in fact formally to a convolution theorem for
Mellin transforms, namely (cf. (1))
$$
{\Cal M}\left[\int_0^\infty f(u)g\bigl({x\over u}\bigr){\d u\over u};s\right]
= \;{\Cal M}[f(x);s]\,{\Cal M}[g(x);s] \;=\; F(s)G(s),\leqno(11)
$$
which is eq. (4.2.22) of I. Sneddon [10]. An alternative proof of (9) follows
from the second formula in (3) and (11), but we need again a result like
Theorem 1 to ensure the validity of the repeated use of (11). A similar approach
via (modified) Mellin transforms and convolutions was carried out by the author
in [7].

\medskip
There is another possibility for the use of the identity (3). Namely, one
can evaluate the Laplace transform of $\{x\}/x$ for real values
of the variable. This is given by

\medskip
THEOREM 3. {\it If $M\ge1$ is a fixed integer and $\g$ denotes Euler's
constant, then for $T\to\infty$}
$$
\eqalign{
\int_0^\infty {\{x\}\over x}{\roman e}^{-x/T}\d x&
= \hf\log T - \hf\g +\hf\log(2\pi)
\cr&
+ \sum_{m=1}^M{\z(1-2m)\over(2m-1)!(1-2m)}T^{1-2m} + O_M(T^{-1-2M}).\cr}\leqno(12)
$$

\medskip
{\bf Proof}. We multiply (3) by $T^s\G(s)$, where $\G(s)$ is the
gamma-function, integrate over $s$ and use the well-known identity
(e.g., see the Appendix of [5])
$$
{\roman e}^{-z} = {1\over2\pi i}\int_{c-i\infty}^{c+i\infty}z^{-s}\G(s)\d s
\qquad(\R z > 0,\,c >0).
$$
We obtain
$$
{1\over2\pi i}\int_{c-i\infty}^{c+i\infty}{\z(s)\over s}\,T^s\G(s)\d s
= -\int_0^\infty {\{x\}\over x}\,{\roman e}^{-x/T}\d x\qquad(0<c<1).\leqno(13)
$$
In the integral on the left-hand side of (13) we shift the line of integration
to $\R s = -N - 1/2, \,N = 2M+1$ (i.e., taking $c = -N - 1/2$) and then apply
the residue theorem. The gamma-function has simple
poles at $s = -m,\,m = 0,1,2,\ldots$ with residues $(-1)^m/m!$. The zeta-function
has simple (so-called ``trivial zeros") at $s = -2m,\,m\in\NN$, which cancel with
the corresponding poles of $\G(s)$. Thus there remains a pole of order two at
$s=0$, plus simple poles at $s = -1,-3,-5,\ldots\,$. The former produces the
main term in (12), when we take into account that
$\z(0) = -\hf,\,\z'(0) = -\hf\log(2\pi)$ (see [5, Chapter 1]) and $\G'(1) = -\g$. The
simple poles at $s = -1,-3,-5,\ldots\,$ produce the sum over $m$ in (12), and the
proof is complete.

\medskip
{\bf Remark 2}. The method of proof clearly yields also, as $T\to\infty$,
$$
\int_0^\infty {\f_n(x)\over x}{\roman e}^{-x/T}\d x
= P_n(\log T) + \sum_{m=1}^M c_{m,n}T^{1-2m} + O_M(T^{-1-2M}),
$$
where $P_n(z)$ is a polynomial in $z$ of degree $n$ whose coefficients may
be explicitly evaluated, and $c_{m,n}$ are suitable constants which also
may be explicitly evaluated.

\medskip
For our last result we turn to M\"untz's identity (4)--(5)  and choose
$f(x) = {\roman e}^{-\pi x^2}$, which is a fast converging kernel
function. Then
$$\eqalign{
Pf(x) &= \sum_{n=1}^\infty f(nx) - \frac{1}{x}\int_0^\infty f(t)\d t
= \sum_{n=1}^\infty {\roman e}^{-\pi n^2x^2} - {1\over2x},\cr
F(s) &= \int_0^\infty {\roman e}^{-\pi x^2}x^{s-1}\d x = \hf \pi^{-s/2}\G(\hf s).\cr}
$$
From (2) and (4) it follows then that, for $0<\s<1$,
$$
\int_{-\infty}^\infty|\z(\s+it)\G(\hf \s+ \hf it)|^2\d t
= 8\pi^{1+\s}\int_0^\infty\left(\sum_{n=1}^\infty {\roman e}^{-\pi n^2x^2} - {1\over2x}
\right)^2x^{2\s-1}\d x.\leqno(14)
$$
The series on the right-hand side of (14) is connected to Jacobi's
theta function
$$
\t(z) := \sum_{n=-\infty}^\infty {\roman e}^{-\pi n^2z}\qquad(\R z > 0),\leqno(15)
$$
which satisfies the functional equation (proved easily by e.g., Poisson
summation formula)
$$
\t(t) = {1\over\sqrt{t}}\,\t\left({1\over t}\right)\qquad(t > 0).\leqno(16)
$$
From (15)--(16) we infer that
$$
\sum_{n=1}^\infty {\roman e}^{-\pi n^2x^2} =
{1\over2}\bigl(\t(x^2) - 1\bigr) = {1\over2x}\t\left({1\over x^2}\right)
- {1\over2}\qquad(x>0). \leqno(17)
$$
By using (17) it is seen that the right-hand side of (14) equals
$$
\eqalign{&
8\pi^{1+\s}\int_0^\infty {1\over4}\left({1\over x}\t\left({1\over x^2}\right)
- 1 - {1\over x}\right)^2x^{2\s-1}\d x\cr&
= 2\pi^{1+\s}\int_0^\infty (u\t(u^2) - 1 - u)^2u^{-1-2\s}\d u.\cr}\leqno(18)
$$
The  (absolute) convergence of the last integral at infinity follows from
$$
u\t(u^2)  - u = 2u\sum_{n=1}^\infty {\roman e}^{-\pi n^2u^2},
$$
while the convergence at zero  follows from
$$
u\t(u^2) = \t\left({1\over u^2}\right) = 1 + O\left({\roman e}^{-u^{-2}}\right)
\qquad(u\to 0+).
$$
Now we note that (14) remains unchanged when $\s$ is replaced by $1-\s$, and then
we use the functional equation (see e.g., [5, Chapter 1]) for $\z(s)$ in the form
$$
\pi^{-s/2}\z(s)\G(\hf s) = \pi^{-(1-s)/2}\z(1-s)\G(\hf(1- s))
$$
to transform the resulting left-hand side of (14). Then (14) and (18) yield
the following
\medskip
THEOREM 4. {\it For $0<\s<1$ we have}
$$
\int_{-\infty}^\infty|\z(\s+it)\G(\hf \s+ \hf it)|^2\d t
= 2\pi^{\s+1}\int_0^\infty(u\t(u^2)-1-u)^2u^{2\s-3}\d u.
$$

\vfill
\eject\topskip2cm

\bigskip\bigskip
\Refs
\bigskip

\item{[1]}
L. B\'{a}ez-Duarte, A strengthening of the Nyman-Beurling criterion for the
Riemann hypothesis, Rendiconti Accad. Lincei, {\bf23} (2003) 5-11.

\item{[2]} L. B\'aez-Duarte,
A general strong Nyman-Beurling criterion for the Riemann Hypothesis,
to appear in Publs. Inst. Math. Belgrade.

\item{[3]} L. B\'{a}ez-Duarte, M. Balazard, B. Landreau et E. Saias,
Notes sur la fonction $\zeta$ de Riemann, 3, Adv. Math. {\bf149} (2000) 130-144.

\item{[4]}A. Beurling, A closure problem related to the Riemann zeta-function,
Proc. Nat. Acad. Sci. U.S.A. {\bf41} (1955), 312-314.

\item{[5]} A. Ivi\'c, The Riemann zeta-function, John Wiley \&
Sons, New York, 1985.

\item{[6]} A. Ivi\'c, Some identities for the Riemann zeta-function,
Univ. Beograd. Publ. Elektrotehn. Fak. Ser. Mat. {\bf14}(2003), 20-25.

\item{[7]} A. Ivi\'c, Estimates of convolutions of certain number-theoretic
error terms, Intern. J. Math. and Mathematical Sciences Vol. 2004, No. 1,
1-23.

\item{[8]} C. H. M\"{u}ntz, Beziehungen der Riemannschen $\zeta$-Funktion zu
willk\"{u}rlichen reellen Funktionen, Mat. Tidsskrift, B (1922), 39-47.

\item{[9]} B. Nyman, On the One-Dimensional Translation Group and
Semi-Group in Certain Function Spaces, Thesis, University of Uppsala, 1950, 55p.

\item{[10]} I. Sneddon, The Use of Integral Transforms,  McGraw-Hill, New York
etc., 1972.

\item{[11]} E.C. Titchmarsh, Introduction to the Theory of Fourier
Integrals, Oxford University Press, Oxford, 1948.

\item{[12]} E.C. Titchmarsh,
The Theory of the Riemann Zeta-Function (2nd ed.),
Oxford at the Clarendon Press, 1986.

\vskip2cm

Aleksandar Ivi\'c \par
Katedra Matematike RGF-a\par
Universiteta u Beogradu\par
\DJ u\v sina 7, 11000 Beograd\par
Serbia and Montenegro\par
{\sevenbf e-mail: aivic\@matf.bg.ac.yu, ivic\@rgf.bg.ac.yu}
\endRefs


\bye